\documentclass[a4pape,11pt]{amsart} 
\usepackage{amsmath,amssymb,amsxtra,amsthm,amscd}
\usepackage[margin=3cm]{geometry} 
\usepackage{hyperref} 
\usepackage{graphicx,color} 
\usepackage{epic,eepic} 
\usepackage[all]{xy}

\newtheorem{theorem}{Theorem}[section]

\theoremstyle{definition}
\newtheorem{example}[theorem]{Example}
 
\newtheorem{remark}[theorem]{Remark}

\newcommand{\C}{\mathbb{C}} 
\newcommand{\R}{\mathbb{R}} 
\newcommand{\PP}{\mathbb{P}} 
\newcommand{\Z}{\mathbb{Z}} 

\newcommand{\hGamma}{\widehat{\Gamma}} 

\newcommand{\iu}{\mathtt{i}}

\newcommand{\cN}{\mathcal{N}} 
\newcommand{\cV}{\mathcal{V}} 
\newcommand{\cL}{\mathcal{L}} 

\newcommand{\id}{\operatorname{id}} 
\newcommand{\ch}{\operatorname{ch}} 
\newcommand{\Coh}{\operatorname{Coh}} 
\newcommand{\Fuk}{\operatorname{Fuk}} 
\newcommand{\FS}{\operatorname{FS}} 
\newcommand{\Log}{\operatorname{Log}} 
\newcommand{\pt}{\operatorname{pt}} 
\newcommand{\Sing}{\operatorname{Sing}} 

\newcommand{\corr}[1]{\left\langle #1 \right\rangle} 
\newcommand{\pair}[2]{\langle #1,#2 \rangle} 
\newcommand{\parfrac}[2]{\frac{\partial #1}{\partial #2}}

\renewcommand{\Im}{\operatorname{Im}} 

\begin{document} 
\title{Gamma conjecture and tropical geometry} 
\date{19 December 2021}
\author{Hiroshi Iritani} 
\email{iritani@math.kyoto-u.ac.jp} 
\address{Department of Mathematics, Graduate School of Science, Kyoto University, 
Kitashirakawa-Oiwake-cho, Sakyo-ku, Kyoto, 606-8502, JAPAN}
\maketitle

\begin{abstract} 
Hodge-theoretic mirror symmetry for a Calabi-Yau mirror pair says that the variation of Hodge structure arising from quantum cohomology of a Calabi-Yau manifold and that arising from deformation of complex structures on the dual Calabi-Yau manifold can be identified with each other, and it has been conjectured ($\hGamma$-conjecture) that the $\hGamma$-integral structure \cite{Iritani:Int} in quantum cohomology corresponds to a natural integral structure on the mirror side. 
Here the $\hGamma$-integral structure is defined via the topological $K$-group and the $\hGamma$-class, a characteristic class with transcendental coefficients containing the Riemann $\zeta$-values. In this article, we explain an approach to the $\hGamma$-conjecture using tropical geometry and observe that the Riemann $\zeta$-values arise as error terms of tropicalization in the computation of mirror periods. This is based on joint work \cite{AGIS} with Abouzaid, Ganatra and Sheridan. 
\end{abstract} 

\section{Mirror symmetry} 
Mirror symmetry is a conjectural duality between symplectic and complex geometry. It roughly speaking predicts that symplectic topology on a symplectic manifold $Y$ ``corresponds'' (or equivalent) to complex geometry on another complex manifold $Z$. Then $Z$ is called a mirror of $Y$ and vice versa.  
In this article, we consider two kinds of mirror correspondences: Calabi-Yau mirror pairs and Fano/LG mirror pairs. 

\begin{description} 
\item[Calabi-Yau mirror pairs] 
\[
(Y, (-\log t) \omega) \quad \longleftrightarrow \quad (Z_t, \Omega_t)_{t\in\Delta^*} 
\]
In the left-hand side, we consider a Calabi-Yau manifold $Y$ equipped with a K\"ahler form $\omega$ whose cohomology class is integral: they give a family $(Y, (-\log t)\omega)$ of symplectic manifolds parametrized by small $t>0$. On the right-hand side we consider a family of Calabi-Yau manifolds $Z_t$ of the same dimension equipped with a holomorphic volume form $\Omega_t$. The family $\{Z_t\}$ maximally degenerates at $t=0$ in a suitable sense: the monodromy $M$ on $H^n_{\rm prim}(Z_t)$ is maximally unipotent (i.e.~$(M-\id)^{n} \neq 0, (M-\id)^{n+1} =0$ for $n=\dim Z_t$) and the limiting mixed Hodge structure is of Hodge-Tate type (see \cite{Deligne:infinity}). The limit point $t=0$ is called the large complex structure limit.  
\item[Fano/LG mirror pairs] 
\[
F \quad \longleftrightarrow \quad W \colon (\C^\times)^n \to \C
\]
In the left-hand side we consider a Fano manifold $F$ (or a monotone symplectic manifold), i.e.~$c_1(F)$ is represented by a positive $(1,1)$-form, which gives a symplectic form. In the right-hand side we consider a Landau-Ginzburg model, a Laurent polynomial function $W$ on the algebraic torus $(\C^\times)^n$ with $n=\dim F$. 
\end{description} 
On the complex geometry side, we consider (exponential) periods, namely, integrals of the following form 
\begin{align*} 
& \int_{C_t \subset Z_t} \Omega_t  \quad && \text{for an $n$-cycle $C_t$ (for  Calabi-Yau mirrors)} \\
& \int_{\Gamma} e^{-t W} \frac{dx_1}{x_1} \cdots \frac{dx_n}{x_n} 
\quad && \text{for a not necessarily compact $n$-cycle $\Gamma$ (for Fano mirrors)}
\end{align*} 
Under mirror symmetry, they should yield solutions to the quantum differential equation on the symplectic side, which is defined by counting rational curves (the genus-zero Gromov-Witten invariants). In the Calabi-Yau case, the mirror correspondence can be formulated as an isomorphism of variation of Hodge structure (VHS) as given in Table \ref{tab:VHS_CYmirrorpairs}. 
\begin{table}[h] 
\caption{Mirror correspondence on the level of VHS (see e.g.~\cite{Iritani:periods}). 
Here $\{\phi_i\}$ is a basis of $H^{1,1}_{\rm alg}(X)$ and $\tau = \sum_{i=1}^r \tau^i \phi_i$.} 
\label{tab:VHS_CYmirrorpairs}
\begin{tabular}{r|c|c}
& symplectic side & complex side \\ 
\hline 
bundle & $\phantom{\Bigg|}$ $\left(\bigoplus_p H^{p,p}_{\rm alg}(Y)\right) \times H^{1,1}_{\rm alg}(Y) \to H^{1,1}_{\rm alg}(Y)$ 
& $\bigcup_t H^n_{\rm prim}(Z_t) \to \Delta^*$ \\ 
connection & quantum connection 
$\nabla = d + \displaystyle\sum_{i=1}^r (\phi_i\star_\tau) d\tau^i$ & 
Gauss-Manin connection $\nabla^{\rm GM}$ \\ 
filtration & $F^p = \displaystyle \bigoplus_{k\le n-p} H^{k,k}_{\rm alg}(Y)$ 
& 
$F^p = \displaystyle \bigoplus_{k\ge p} H^{n-k,k}_{\rm prim}(Z_t)$ \\ 
polarization & $\displaystyle (2\pi\iu)^n \int_Y ((-1)^{\deg/2} \alpha) \cup \beta$ 
& $\displaystyle (-1)^{\frac{n(n-1)}{2}} \int_X \alpha \cup \beta$ \\ 
$\Z$-structure & $\phantom{\Big|}$ $\hGamma$-integral structure 
$K^0_{\rm alg}(X) \to \bigoplus_{k} H^{k,k}_{\rm alg}(X)$, 
& $H^n(Z_t,\Z)$ \\
& $V \mapsto \hGamma_X (2\pi\iu)^{\deg/2} \ch(V)$ & 
\end{tabular}  
\end{table} 

\begin{remark} 
In Table \ref{tab:VHS_CYmirrorpairs}, we restrict our attention to the algebraic part of cohomology in the symplectic side, and the primitive part of the middle cohomology in the complex side. The $\hGamma$-integral structure is also restricted to the $K$-group of algebraic vector bundles. We do not know much about Hodge-theoretic mirror symmetry beyond these parts. On the symplectic side, we can naturally define the $\hGamma$-integral structure corresponding to the topolgoical $K$-group \cite{Iritani:Int}, so it might be possible to consider an integral structure on the complex side corresponding to the topological $K$-group as well. A conjectural isomorphism $K_{\rm top}^*(Y) \cong K_{\rm top}^*(Z)$ has been discussed in the literature \cite{Treuman}. 
\end{remark} 

\section{The $\hGamma$-class} 
Let $X$ be an almost complex manifold. We write the total Chern class of the tangent bundle $TX$ as 
\[
c(TX) = (1+\delta_1) (1+\delta_2) \cdots (1+\delta_n)
\]
where $\delta_1,\dots,\delta_n$ are virtual cohomology classes called the \emph{Chern roots}. Each $\delta_i$ may not exist as a cohomology class, but any symmetric functions in $\delta_1,\dots,\delta_n$ can be written as polynomials in $c_1(TX),\dots,c_n(TX)$ and make sense as cohomology classes.  
The \emph{$\hGamma$-class} of $X$ is defined to be 
\[
\hGamma_X = \Gamma(1+\delta_1) \cdots \Gamma(1+\delta_n) \in H^*(X,\R) 
\]
where $\Gamma(1+x) = \int_0^\infty e^{-t} t^x dt/t$ is the Euler $\Gamma$-function. By the Taylor expansion of the $\hGamma$-function, we can think of the right-hand side as a symmetric power series of $\delta_1,\dots,\delta_n$: then the right-hand side makes sense as a cohomology class of $X$. The $\hGamma$-function $\Gamma(1+x)$ has simple poles at $x=-1,-2,-3,\dots$ and it has the following infinite product expansion:  
\[
\Gamma(1+x) = \frac{e^{-\gamma x}}{\prod_{n\ge 0} (1+x/n) e^{-x/n}}   
\]
where $\gamma = \lim_{n\to \infty} (1+ \frac{1}{2} + \cdots + \frac{1}{n} - \log n)$. This can be calculated as 
\begin{align*} 
\Gamma(1+x) & = e^{-\gamma x} \prod_{n=1}^\infty e^{-\log (1+ \frac{x}{n}) + \frac{x}{n}} \\
& = e^{-\gamma x} \prod_{n=1}^\infty \exp\left( \sum_{k=2}^\infty \frac{(-1)^k}{k} \frac{x^k}{n^k}\right) \\ 
& = \exp\left(-\gamma x + \sum_{k=2}^\infty (-1)^k \frac{\zeta(k)}{k}  x^k \right) 
\end{align*} 
where $\zeta(k) = \sum_{n=1}^\infty \frac{1}{n^k}$ is the value of the Riemann $\zeta$-function. Hence the $\hGamma$-class can be written as 
\[
\hGamma_X=  \exp\left( -\gamma c_1(X) + \sum_{k=2}^\infty (-1)^k 
\frac{\zeta(k)}{(k-1)!} \ch_k(TX) \right) 
\]
\begin{remark} 
The $\hGamma$-class has the following geometric interpretation. We consider a free loop space $LX$ equipped with the $S^1$-action rotating loops. We also consider the set $X$ of constant loops in $LX$. Then the $\hGamma$-class can be interpreted as a regularization of the $S^1$-equivariant Euler class of the positive normal bundle $\cN_+$ of $X$ in $LX$ \cite{Lu,GGI}; here positive means the positive weight part as an $S^1$-representation. We have 
\[
(2\pi)^{n/2} z^{n-\frac{\deg}{2}} z^{c_1(X)} \hGamma_X \sim 
e_{S^1}(\cN_+) = \frac{1}{\prod_i \prod_{k \ge 0} (\delta_i + k z)}. 
\]
\end{remark} 

\section{Mirror symmetric $\hGamma$-conjecture} 
The mirror symmetric  $\hGamma$-conjecture roughly speaking says that the $\hGamma$-integral structure \cite{Iritani:Int} on the symplectic side should correspond to a natural integral structure on the complex side. In this article, we do not discuss Hodge-theoretic mirror symmetry with integral structure anymore (we do not even give the definition of quantum cohomology or the $\hGamma$-integral structure on it): we refer the reader to e.g.~\cite{Iritani:sugaku, GI, Iritani:periods}. Instead, we discuss a more concrete conjecture, ``mirror symmetric $\hGamma$-conjecture'' stated in terms of (exponential) periods. This problem is more of a topological nature, and does not involve counting rational curves. 

The mirror symmetric $\hGamma$-conjecture originates from Hosono's conjecture \cite{Hosono:centralcharges}, which says that mirror periods equal the pairing of certain explicit hypergeometric series with the Chern classes of vector bundles, for Batyrev mirror pairs of Calabi-Yau hypersurfaces. By taking the asymptotics of Hosono's conjecture at the large complex structure limit, we arrive at the following conjecture. 
\begin{description} 
\item[Calabi-Yau case] Let $(Y,(-\log t) \omega)$ and $(Z_t, \Omega_t)$ be a Calabi-Yau mirror pair as above. For a certain family of $n$-cycles $C_t \subset Z_t$, we have a $K$-theory class $V$ on $Y$ such that 
\begin{equation} 
\label{eq:conj_CY} 
\int_{C_t\subset Z_t} \Omega_t = \int_Y t^{-\omega} \hGamma_Y  (2\pi\iu)^{\deg/2} \ch(V) + O(t^\epsilon) 
\end{equation} 
as $t\to +0$, where $\epsilon>0$ is a positive real number. 

\item[Fano/LG case] Let $F$ and $W\colon (\C^\times)^n \to \C$ be a Fano/LG mirror pair as above. For a certain (possibly noncompact) $n$-cycle $\Gamma \subset (\C^\times)^n$, we have a $K$-theory class $V$ on $F$ such that 
\begin{equation} 
\label{eq:conj_Fano}
\int_\Gamma e^{-tW} \frac{dx_1\cdots dx_n}{x_1\cdots x_n} 
= \int_F t^{-c_1(F)} \hGamma_F (2\pi\iu)^{\deg/2} \ch(V) + O(t^\epsilon) 
\end{equation} 
as $t\to +0$, where $\epsilon>0$ is a positive real number. 
\end{description} 

This conjecture has been verified for (weak) Fano toric orbifolds and certain complete intersections in them, for some choices of $V$ and cycles $C_t$ (or $\Gamma$), see  \cite{Iritani:Int, Iritani:periods}. The paper \cite{AGIS} gives another proof for Batyrev mirror pairs based on the SYZ picture and tropical geometry. 

\begin{remark} 
(A) 
This conjecture is closely related to homological mirror symmetry (with symplectic side and complex side interchanged). Homological mirror symmetry predicts that the derived category of coherent sheaf on one side should be equivalent to the Fukaya (or Fukaya-Seidel) category of the other side: 
\begin{align*} 
D^b\Coh (Y) & \cong D^b \Fuk(Z_t)  && \text{for Calabi-Yau mirror pairs} \\
D^b\Coh (F) & \cong D^b \FS(W)  && \text{for Fano/LG mirror pairs} 
\end{align*} 
We expect that the $K$-classes $V$ and the cycles $C_t$ (or $\Gamma$) in the above conjecture should correspond to each other under homological mirror symmetry. Namely, when $V$ comes from a coherent sheaf $\cV\in D^b\Coh(Y)$ or $D^b\Coh(F)$, the cycle $C_t$ (or $\Gamma$) should be the Lagrangian submanifold $\cL$  mirror to $\cV$. 

(B) The categorical equivalence in homological mirror symmetry is given only up to auto-equivalences, and hence the correspondence in (A) is ambiguous. In this conjecture, more precisely, we should consider the equivalence induced from Strominger-Yau-Zaslow (SYZ) dual torus fibrations (see \cite{SYZ,Gross:SYZ,Gross:SYZ2}). The SYZ conjecture says that (in the Calabi-Yau case) we have special Lagrangian torus fibrations\footnote{For tropical computation of periods, we do not need Ricci-flat metrics or special Lagrangian fibrations: we only need a weaker version as in the Gross-Siebert program  \cite{Gross-Siebert}.} $p_1\colon Y \to B$, $p_2 \colon Z_t \to B$ with singularities 
\[
\xymatrix{
Y \ar[dr]^{p_1} & & Z_t \ar[dl]_{p_2} \\
& B & 
}
\]
that are dual to each other, where $B$ is a real $n$-dimensional manifold homeomorphic to a sphere. Here $Y$ and $Z_t$ are equipped with Ricci-flat K\"ahler metrics. It is expected that $Z_t$ converges\footnote{Likewise, we can consider a maximal degeneration of complex strcutures on $Y$, which corresponds under mirror symmetry to the large-radius limit for the K\"ahler (symplectic) structure on $Z_t$, and the complex degeneration induces the collapse $Y\to B$. When we take into account both the symplectic and complex structures, we should consider a mirror pair $(Y_s, (-\log t) \omega_Y) \leftrightarrow (Z_t, (-\log s) \omega_Z)$ of maximally degenerating families.}, as $t\to 0$, to the base $B$ in the sense of Gromov-Hausdorff topology (where we normalize the metric so that the diameter is constant). The SYZ fibrations should induce  categorical equivalences as above\footnote{More precisely, the categorical equivalence would be given up to the twist by a line bundle, and the ambiguity would be fixed by choosing a Lagrangian section that corresponds to the structure sheaf.}, and it is expected that a Lagrangian section of $p_2 \colon Z_t \to B$ corresponds to a line bundle on $Y$. In particular, if we take a Lagrangian section $C_{0,t}$ (or $\Gamma_0$) corresponding to the structure sheaf, we have 
\begin{align} 
\label{eq:gammaconj}
\begin{split} 
\int_{C_{0,t}\subset Z_t} \Omega_t & 
= \int_Y t^{-\omega} \hGamma_Y  + O(t^\epsilon) \qquad \text{(in the Calabi-Yau case)} \\ 
\int_{\Gamma_0} e^{-tW} \frac{dx_1\cdots dx_n}{x_1\cdots x_n} 
& = \int_F t^{-c_1(F)} \hGamma_F + O(t^\epsilon)\qquad \text{(in the Fano case)} 
\end{split} 
\end{align} 
as $t\to +0$. 
In examples such as Batyrev mirrors, $C_{0,t}$ arises as a ``positive real locus''. In the Fano case, we expect that $\Gamma_0 = (\R_{>0})^n$. 

When $V$ is the class of the structure sheaf of a point, the corresponding cycle should be a fibre $p_2^{-1}(b)$ of the SYZ fibration $p_2\colon Z_t\to B$ in the Calabi-Yau case, and the compact torus $(S^1)^n\subset (\C^\times)^n$ in the Fano case. Since $\ch(V) = [\pt]$ in this case, we do not see nontrivial components of the $\hGamma$-class in the asymptotics of the corresponding (exponential) periods. Tonkonog \cite{Tonkonog} showed that the exponential period $\int_{(S^1)^n} e^{-t W} \frac{dx_1\cdots dx_n}{x_1\cdots x_n}$ of $(S^1)^n$ is a generating series of gravitational descendants of $F$ (the pairing of the $J$-function and $[\pt]$), when one chooses a monotone Lagrangian submanifold $L$ in $F$ and defines $W$ by counting holomorphic discs with boundaries in $L$. See also (C) below. 

(C) The above conjectures \eqref{eq:conj_CY}, \eqref{eq:conj_Fano}, \eqref{eq:gammaconj} say that the $\hGamma$-class appears in the asymptotics of (exponential) periods of the mirror in the large complex structure limit $t\to 0$. The leading asymptotics are polynomials in $\log t$. The information of curve counting is contained in the higher order terms in $t$, which are exponentially small compared to the asymptotic part. If we include all the higher-order terms, the right-hand side in the conjecture should become 
\begin{align*} 
& \int_Y J_Y( \omega \log t,-1) \cup \hGamma_Y (2\pi\iu)^{\deg/2} \ch(V) && \text{(in the Calabi-Yau case)}\\ 
& \int_F J_F( c_1(F) \log t, -1) \cup \hGamma_F (2\pi\iu)^{\deg/2} \ch(V) && 
\text{(in the Fano case)}
\end{align*} 
where $J_X(\tau,z)$ is the small $J$-function 
\[
J_X(\tau,z) = e^{\tau/z} \left( 1 + \sum_i \sum_{d\in H_2(X,\Z)} 
\corr{\frac{\phi^i}{z(z-\psi)}} e^{\pair{\tau}{d}} \phi_i \right)  
\]
defined in terms of gravitational Gromov-Witten invariants (see \cite{GGI,GI,Iritani:sugaku} for the notation). 
In the Calabi-Yau case, we need to normalize the volume form $\Omega_t$ by a  Hodge-theoretic condition, as discussed in \cite{Deligne:infinity}; we also need to assume that the parameter $t$ of the mirror family $\{Z_t\}$ is normalized so that the mirror map is trivial $\tau= \omega \log t$. 
\end{remark} 

\begin{example} 
Let $Y$ be a Calabi-Yau 3-fold. The asymptotic part of the first equation of \eqref{eq:gammaconj} takes the form: 
\[
\int_Y t^{-\omega} \hGamma = (-\log t)^3 \underbrace{\int_Y \frac{\omega^3}{3!}}_{\text{volume of $Y$}}   
- (-\log t) \zeta(2) \int_Y \omega \cup c_2(Y) 
- \zeta(3) \underbrace{\int_Y c_3(Y)}_{\text{Euler number}}  
\]
Note that the leading term is the symplectic volume of $(Y, (-\log t) \omega)$. 
\end{example} 

\begin{example} 
Let $F$ be the projective space $\PP^n$ and let $Y$ be a degree $n+1$ Calabi-Yau hypersurface in $F=\PP^n$. The mirror of $F$ is given by the Laurent polynomial 
\[
W = x_1+\cdots+x_n+ \frac{1}{x_1\cdots x_n}
\] 
and the mirror of $Y$ is given by a Calabi-Yau compactification $Z_t$ of the fibre $W^{-1}(1/t)$ equipped with the holomorphic volume form 
\[
\Omega_t =\left. \frac{\frac{dx_1}{x_1} \wedge \cdots \wedge \frac{dx_n}{x_n}}{t\cdot dW} \right|_{W^{-1}(1/t)}. 
\]
The mirror symmetric $\hGamma$-conjecture for these mirror pairs holds when $V$ is the structure sheaf and the cycle is the positive real locus \cite{Iritani:Int,Iritani:periods}. We have 
\begin{align*} 
\int_{Z_t \cap (\R_{>0})^n} \Omega_t & = \int_Y t^{-c_1(F)} \hGamma_Y +O(t^\epsilon) \\ 
\int_{(\R_{>0})^n} e^{-t W} \frac{dx_1\cdots dx_n}{x_1\cdots dx_n} & =  
\int_F t^{-c_1(F)} \hGamma_F +O(t^\epsilon) 
\end{align*}
for any $\epsilon$ with $0<\epsilon<n+1$. The higher order term can be also given explicitly as hypergeometric series. We remark that the exponential period in the second line can be written as a Fourier transform of the period in the first line. (The first identity was obtained from the second one by inverse Fourier transformation in \cite{Iritani:periods}.) 
\end{example}

\begin{remark} 
The $\hGamma$-conjecture for Fano manifolds \cite{GGI} proposed by Galkin, Golyshev and the author  does not rely on mirror symmetry (although it sometimes follows from the mirror symmetric $\hGamma$-conjecture). It is formulated purely in terms of quantum cohomology of a Fano manifold. The $\hGamma$-conjecture in \cite{GGI} is related to the Stokes structure of the quantum connection. 
\end{remark}

\section{Periods via tropical geometry} 
We explain how to compute the asymptotics \eqref{eq:gammaconj} of periods using tropical geometry.  Firstly we see that the term 
\[
(-\log t)^n \int_Y \frac{\omega^n}{n!} 
\]
should appear as the leading term of the period $\int_{C_{0,t}} \Omega_t$ using the SYZ picture. Away from the descriminant locus (the singular locus of the fibration), the base $B$ is equipped with a $\Z$-affine structure (i.e.~an atlas as a topological manifold such that every coordinate change is $\Z$-affine linear, that is, belongs to $GL(n;\Z) \ltimes \R^n$) and the SYZ dual fibrations are locally modelled on 
\[
\xymatrix{
T^*B/\Lambda^*  \ar[dr]^{p_1} & & TB/\Lambda \ar[dl]_{p_2} \\ 
 & B & 
}
\]
where $\Lambda \subset TB$ is the lattice defined by the $\Z$-affine structure. 
Let $x_1,\dots,x_n$ be $\Z$-affine local coordinates on $B$. Then the symplectic form on $T^*B/\Lambda^*$ is given by $(- \log t) \omega = \sum_{j=1}^n dx_j \wedge dy^*_j$ and the holomorphic volume form on $TB/\Lambda$ is given by $\Omega=dz_1 \wedge \cdots \wedge dz_n$ with $z_j = x_j + \iu y_j$, where $y^*_j$ are fibre coordinates on $T^*B$ dual to $dx_j$ and $y_j$ are fibre coordinates on $TB$ dual to $\partial/\partial x_j$. Hence we have 
\[
p_{1*} \left(\frac{((-\log t) \omega)^n}{n!} \right) = \text{affine volume form $dx_1\wedge \cdots \wedge dx_n$ on $B$} = s^* \Omega 
\]
where $s$ is a section of $p_2$. From this we expect the leading asymptotics: 
\[
\int_{C_{0,t}} \Omega_t \sim  (-\log t)^n \int_Y \frac{\omega^n}{n!} \qquad 
\text{as $t \to +0$}.   
\]
The error term of this approximation arises from the discriminant locus of the SYZ fibration. 

We explain that tropicalization gives an approximate SYZ(-like) fibration and how we can compute periods tropically by means of example. 

\begin{figure}[ht] 
\begin{minipage}[b]{0.45\hsize} 
\centering 
\includegraphics[scale=0.5,bb=0 20 200 250]{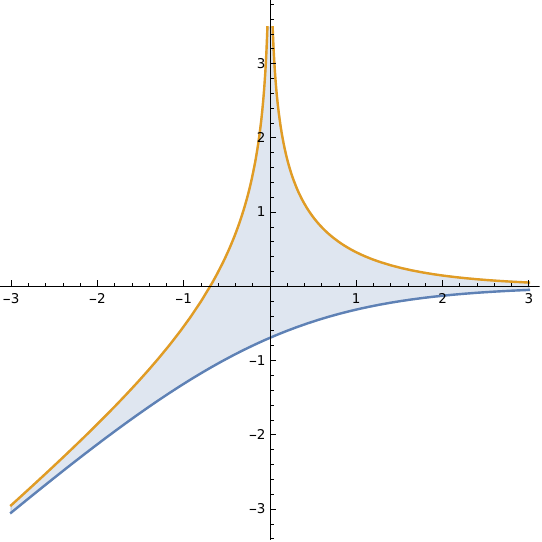}
\caption{Amoeba}
\label{fig:amoeba}  
\end{minipage} 
\begin{minipage}[b]{0.45\hsize}
\centering 
\includegraphics[bb=200 600 400 700]{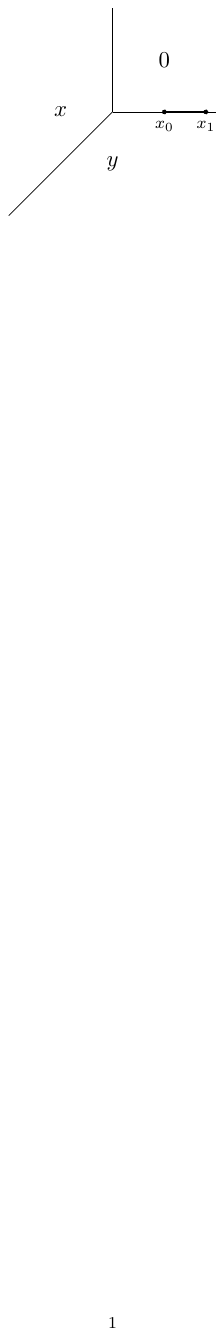} 
\caption{Tropical amoeba} 
\label{fig:tropical_amoeba} 
\end{minipage} 
\centering 
\includegraphics[scale=0.5,bb=150 0 400 300]{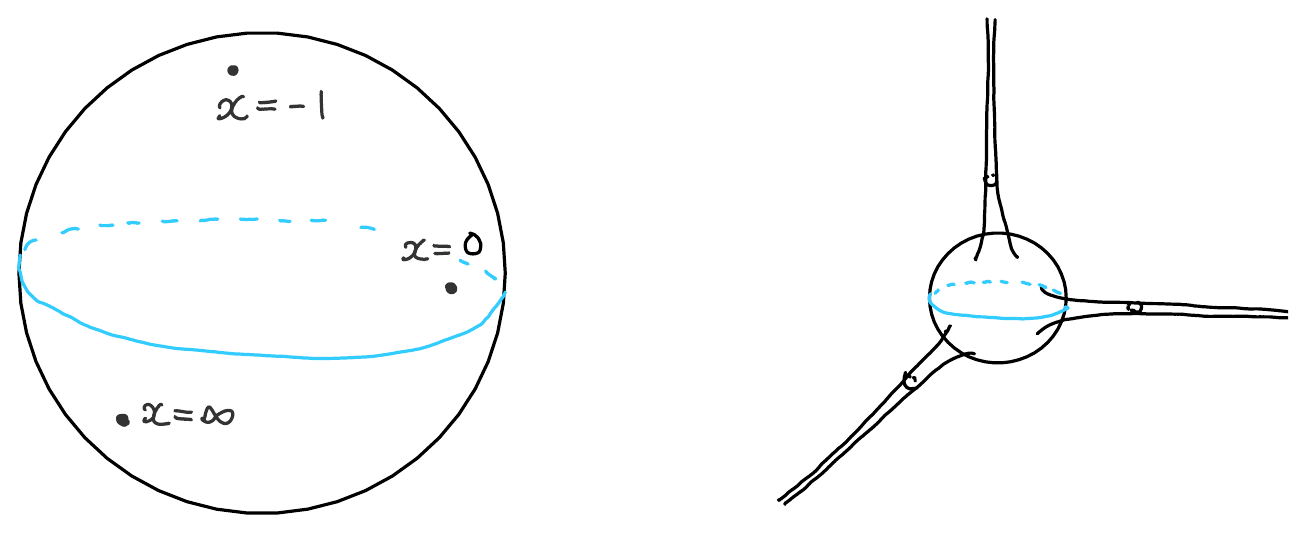}
\caption{Pair of pants ($\PP^1$ minus three points, on the left) and stretching the neck (on the right)}
\label{fig:pairofpants}
\end{figure} 

\begin{example} 
We start with a simple example of tropicalization. Consider the variety 
\[
P = \{(X,Y) \in (\C^\times)^2 : X+Y+1=0\}.  
\]  
This is $\PP^1$ minus three points, and is called a `pair of pants'. We consider the family of maps for $t>0$: 
\[
\Log_t \colon P \to \R^2, \quad (X,Y) \mapsto (x,y) = (\log_t|x|, \log_t|y|). 
\]
The image of the map is called amoeba. In this example, the image is determined by the triangle inequality $||X|-1|\le |Y| = |X+1| \le |X|+1 \Leftrightarrow |t^x-1| \le t^y \le |t^x+1| \Leftrightarrow \log_t(1+t^x)\le y \le \log_t|1-t^x|$. See Figure \ref{fig:amoeba}. If we take the limit $t\to 0$, this approaches to the tropical amoeba (tropical variety) in Figure \ref{fig:tropical_amoeba}. It is given by the singular locus of the piecewise linear function $\min(x,y,0)$. We can think of $\Log_t$ as an approximate torus fibration of $P$ over the tropical amoeba. For example, the fibre at $(x_0,0)$ (with $x_0>0$) is approximately the circle $|X| = t^{x_0}$, $Y = -1$. Most of the points on $P$ go very close to the origin, and three `neck' neighbourhoods of the missing points $X=0$, $X=-1$, $X=\infty$ are stretched to semi-infinite cylinders (see Figure \ref{fig:pairofpants}). 
The integral of a holomorphic 1-form 
\[
\Omega = \left. \frac{\frac{dX}{X}\wedge \frac{dY}{Y}}{d (X+Y+1)}\right|_{X+Y+1=0}
\] 
on a section over the interval $[x_0,x_1]$ in the $x$-axis is approximately the length of the interval multiplied by $(-\log t)$. 
\[
\int_{[x_0,x_1]} s^*\Omega \approx  \int_{[x_0,x_1]} \left. \frac{dX}{XY}\right|_{X=t^x,Y=-1} = (- \log t) (x_1-x_0). 
\]
\end{example} 

\begin{example} 
\label{exa:elliptic} 
Next we consider a tropicalization of an elliptic curve. Let $E_t$ be the affine elliptic curve: 
\[
E_t = \{(X,Y)\in (\C^\times)^2: t(X+Y+\frac{1}{XY}) = 1\} 
\]
This is a compact elliptic curve minus 3 points. The limit of $\Log_t(E_t)$ as $t\to 0$ (tropical amoeba) and the approximate torus fibration are depicted in Figure \ref{fig:tropical_elliptic}. The tropical amoeba is the singular locus of $\min(0,x+1,y+1,1-x-y)$. We see from this picture that $E$ is composed of three `pairs of pants' in the previous example. This is an instance of pairs-of-pants decomposition by Mikhalkin \cite{Mikhalkin:pants}. The approximate fibration given by $\Log_t$ is ``singular'' at three vertices, but we can make it a smooth fibration by adding three missing points (at infinity) to $E_t$ and contracting the unbounded edges by a linear projection around each vertex. For example, around the vertex $(-1,-1)$ of the tropical amoeba, we can use the projection $x-y=\log_t |X| - \log_t |Y|$ along the ray $\R_{\ge 0} (-1,-1)$ to define the fibration. Then we get a smooth $S^1$-fibration over the boundary $\partial \Delta^2\cong S^1$ of the 2-simplex (shown in blue colour). 

\begin{figure}[ht]  
\centering 
\includegraphics[scale=0.7,bb=220 660 300 820]{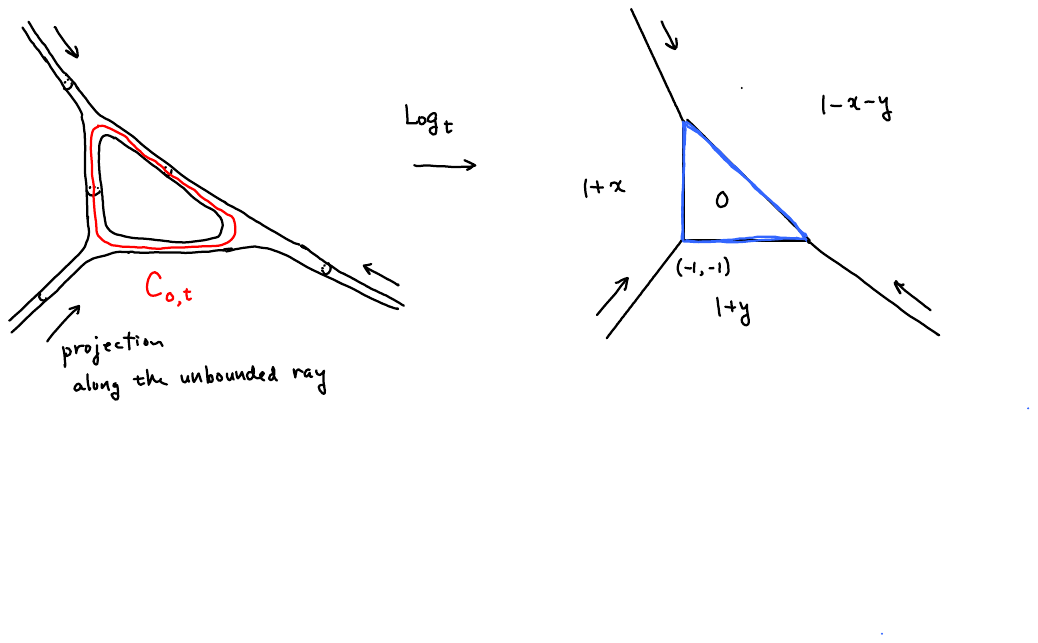}
\caption{Tropical elliptic curve}
\label{fig:tropical_elliptic} 
\end{figure}

The integral of the holomorphic 1-form 
\[
\Omega_t =\left. \frac{d\log X \wedge d\log Y}{d(t(X+Y+1/(XY)))}\right |_{E_t}
\]
over the compact cycle $C_{0,t}$ (shown in Figure \ref{fig:tropical_elliptic}) is asymptotically the same as the affine length ($=9 = 3\times 3$) of the cycle in the tropical base (shown in blue colour) multiplied by $(-\log t)$. 
This elliptic curve $E_t$ is mirror to a cubic curve $Y$ in $\PP^2$. A precise choice of the cubic curve $Y$ does not matter, but $Y$ can be for example:  
\[
Y_s = \{[z_0,z_1,z_2] \in \PP^2 : s(z_0^3 + z_1^3 + z_2^3) + z_0z_1z_2 = 0\}.   
\]
The dual SYZ fibration for $Y_s$ is approximately given by the restriction of the moment mapping $\mu \colon \PP^2 \to \R^2$ (with respect to the anticanonical class $-K_{\PP^2}$ on $\PP^2$) to $Y_s$. 
\[
\mu([z_0,z_1,z_2]) = \left( \frac{3|z_1|^2}{|z_0|^2+|z_1|^2+|z_2|^2}, 
\frac{3 |z_2|^2}{|z_0|^2 + |z_1|^2 + |z_2|^2}\right) 
\]
Note that the image of $\mu$ is the 2-simplex $\{(x,y):x\ge 0, y\ge 0, x+y\le 3\}$ and $\mu(Y_s)$ approaches its boundary as $s\to 0$. The affine length of the tropical cycle equals the symplectic volume $\int_{Y_s} c_1(\PP^2) = 9$. 
\end{example} 

\begin{example} 
\label{exa:2d_local} 
Now we consider the case where the SYZ fibration has singularities. In dimension two, the following local model of singularities appears in the Gross-Siebert program (see \cite{Gross:I,Kontsevich-Soibelman}). 
\begin{align*} 
Z & = \{(X_1,X_2,Y)\in \C^2\times \C^\times : X_1X_2 =1+ Y \} = \C^2_x \setminus \{X_1X_2=1\}
= \PP^2 \setminus (\text{line} \cup \text{conic})\\ 
\Omega & = \frac{d\log X_1 \wedge d \log X_2 \wedge d \log Y}{d\left(\frac{1+Y}{X_1X_2}\right)} 
= \frac{dX_1dX_2}{1-X_1X_2} = \frac{dX_1}{X_1} \wedge \frac{dY}{Y} = \frac{dX_2}{X_2} \wedge \frac{dY}{Y} \\ 
\omega & = \frac{\iu}{2} (dX_1\wedge d\overline{X_1} + dX_2 \wedge d\overline{X_2})
\end{align*} 
This is a log CY variety, in the sense that $\Omega$ has log poles along the boundary divisor $\PP^2\setminus Z = \text{line} \cup \text{conic}$. The Lagrangian\footnote{This is not special Lagrangian since $\Omega \wedge \overline{\Omega}$ is not a constant multiple of $\omega^2$. But we have $\omega|_{T_{r,\lambda}}=0$ and $\Im \Omega|_{T_{r,\lambda}}=0$. Therefore the base of the torus fibration has a $\Z$-affine structure (complex affine structure) defined by fluxes of $\Im \Omega$.} torus fibration on $Z$ is given as follows. We first consider the sympletic reduction by the $S^1$-action $(X_1,X_2,Y) \mapsto (e^{\iu\theta} X_1, e^{-\iu\theta} X_2, Y)$ at the level $\lambda$. The symplectic reduction is identified with the $Y$-plane $\C^\times$. 
\[
\{(X_1,X_2,Y) \in Z : |X_1|^2-|X_2|^2 = \lambda\}/S^1\cong \{Y\in \C^\times\}
\]
By pulling back the standard Lagrangian torus fibration $\{(|Y|=r)\}_{r>0}$ on $\C^\times$, we see that $Z$ is foliated by the family of Lagrangian tori 
\[
T_{\lambda,r} = \{(X_1,X_2,Y) \in Z : |X_1|^2 - |X_2|^2 = \lambda, |Y|=r\} 
\]
parametrized by $\lambda \in \R$ and $r>0$. This gives a torus fibration on $Z$. 
Here $T_{0,1}$ is a unique singular fibre (pinched torus) with singularity at $(0,0,-1) \in T_{0,1}$. 

We compare this fibration with the tropicalization map 
\[
\Log_t \colon Z \to \R^3, \quad (X_1,X_2,Y) \mapsto (\log_t|X_1|, \log_t|X_2|, \log_t|Y|)
\]
for $0<t\ll 1$. 
The last coordinate $\log_t|Y|$ is constant on $T_{\lambda,r}$, but the first two $\log_t|X_1|$, $\log_t|X_2|$ are not. We see however that the map $\Log_t$ approximates the torus fibration away from the singularity $X_1=X_2=0$, i.e.~away from the region where both $|X_1|$ and $|X_2|$ are small. 
When we set $x_1 = \log_t|X_1|$, $x_2=\log_t|X_2|$, $y= \log_t|Y|$, for positive $\epsilon>0$, 
\begin{align*}
\text{$|y|>\epsilon$} \quad & \Longrightarrow \quad \text{$x_i\approx\textstyle  \frac{1}{2}\log_t\left(\frac{\pm \lambda + \sqrt{\lambda^2+4r^2}}{2}\right)$ on $T_{\lambda,r}$} \\ 
\text{$x_1<x_2-\epsilon$} \quad & \Longrightarrow \quad \text{$|X_1|\gg |X_2|$} 
\quad \Longrightarrow  \quad \text{$x_1 \approx \textstyle \frac{1}{2} \log_t |\lambda|$ on $T_{\lambda,r}$} \\ 
\text{$x_2<x_1-\epsilon$} \quad & \Longrightarrow \quad \text{$|X_2|\gg |X_1|$} \quad \Longrightarrow \quad \text{$x_2\approx \textstyle \frac{1}{2}\log_t |\lambda|$ on $T_{\lambda,r}$.} 
\end{align*} 
The image $\Log_t(Z)$ is approximated by the tropical amoeba, which is the singular locus of $\min(0,x_1+x_2,y)$, as in Figure \ref{fig:2d_sing}. From this picture we can see that the above three regions $(|y|>\epsilon)$, $(x_1<x_2-\epsilon)$, $(x_2<x_1-\epsilon)$ cover the region away from the singularity or the origin  $(x_1,x_2,y)=(0,0,0)$. Therefore  
\begin{align*} 
&\text{the map $(x_1,y)$ approximates the torus fibration in the region $(|y|>\epsilon) \cup (x_1<x_2-\epsilon)$;} \\ 
& \text{the map $(x_2,y)$ approximates the torus fibration in the region $(|y|>\epsilon) \cup (x_2<x_1-\epsilon)$.} 
\end{align*} 
Moreover these maps define (approximately) affine-linear coordinates in the respective regions. The map $(x_i,y)$ contracts the face $(x_1+x_2\ge 0=y)$ to the line $y=0$ and send the union $(x_1+x_2=\min(0,y))$ of the other two faces isomorphically to $\R^2$. These coordinates $(x_1,y)$, $(x_2,y)$ are glued together as in Figure \ref{fig:focus-focus}: the base of the Lagrangian fibration has the focus-focus singularity with monodromy $\scriptsize \begin{pmatrix} 1 & 1 \\ 0 & 1 \end{pmatrix}$ (see \cite{Kontsevich-Soibelman,Yamamoto}). 

\begin{figure}[ht] 
\centering
\includegraphics[scale=0.5, bb=80 0 600 434]{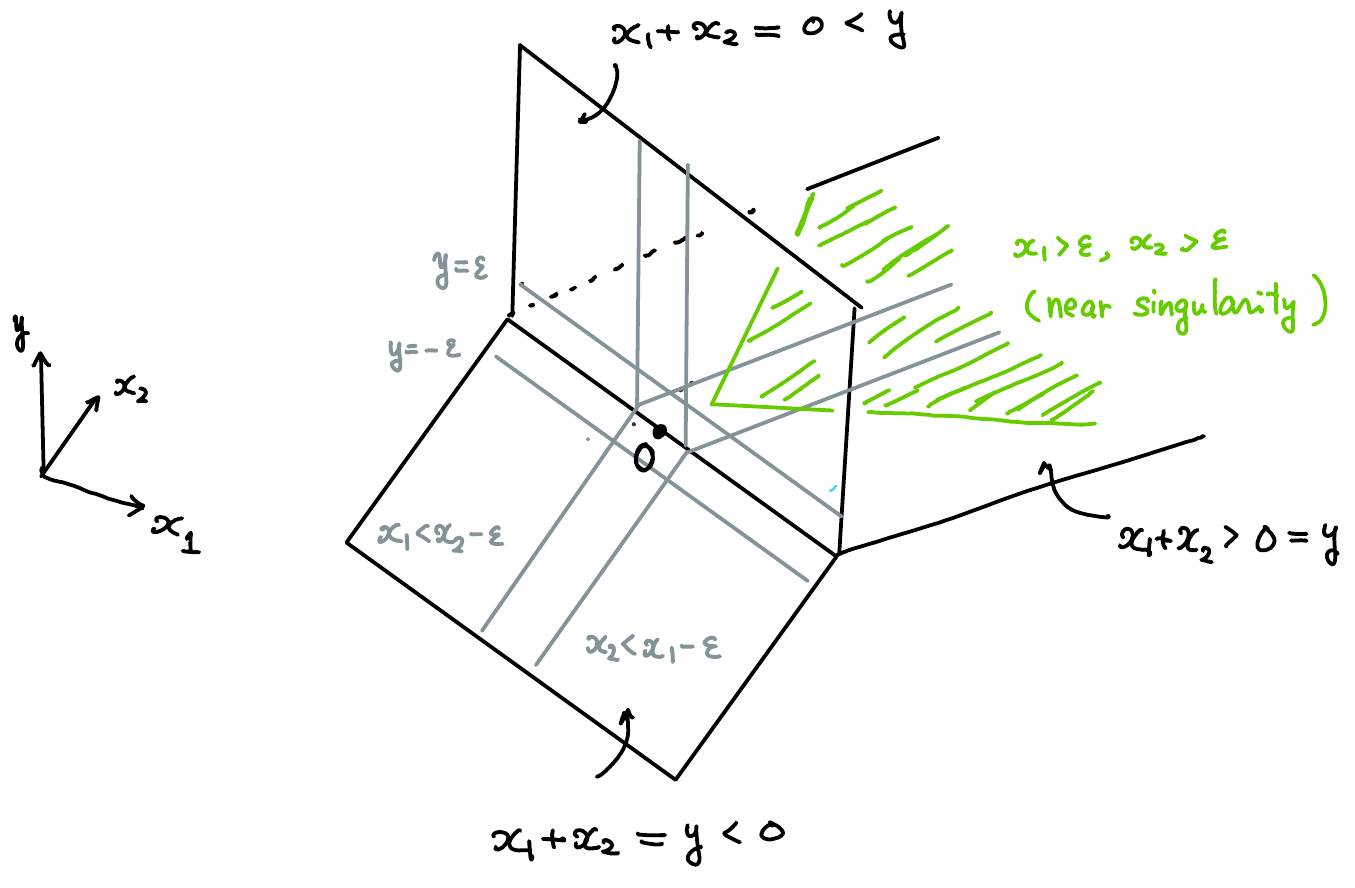} 
\caption{Tropical amoeba of $Z$}
\label{fig:2d_sing}
\end{figure}

\begin{figure}[ht]
\centering
\includegraphics[bb=200 600 400 700]{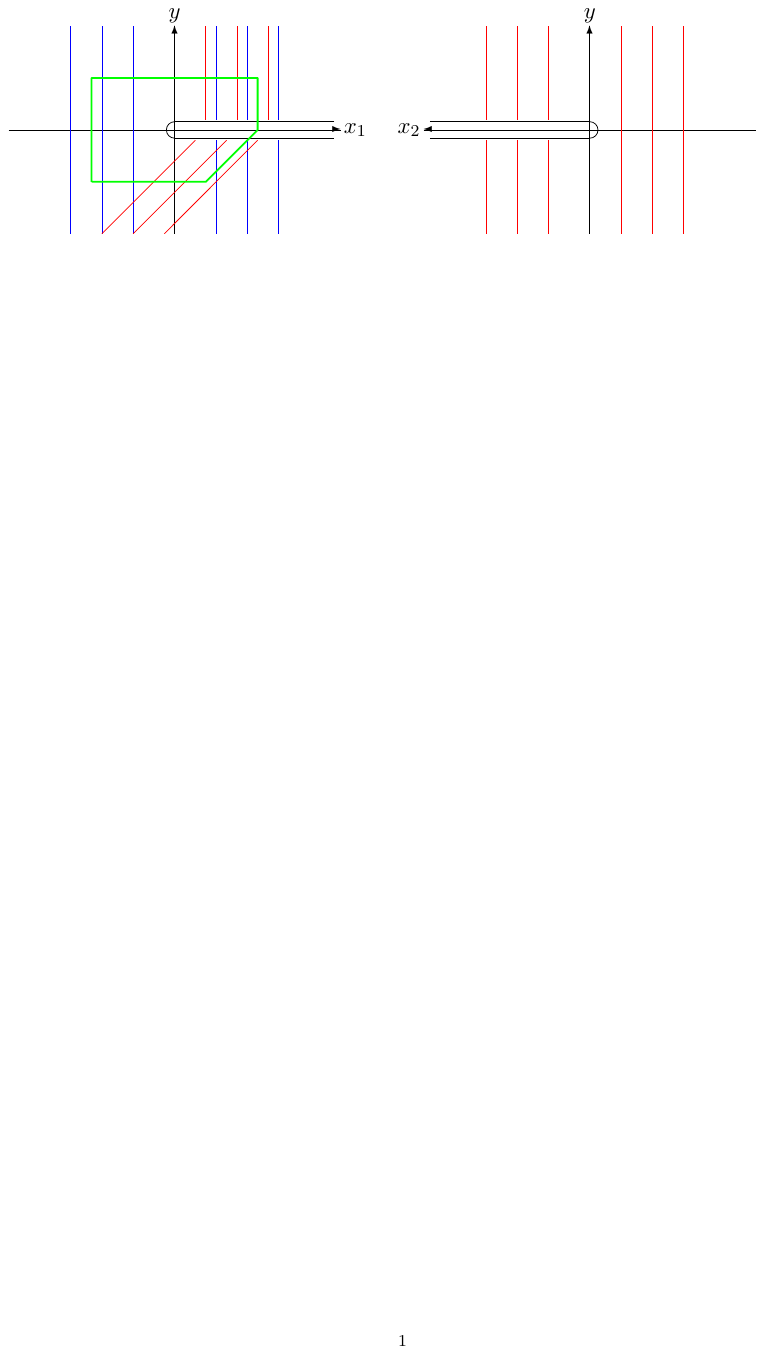}
\caption{Two approximately affine-linear charts $(x_1,y)$ and $(x_2,y)$. Each chart $(x_i,y)$ is affine-linear away from the positive $x_i$-axis. The coordinate change between the charts is given by $x_1+x_2 = \min(0,y)$. The green line shows the polytope corresponding to $R$.}
\label{fig:focus-focus} 
\end{figure} 

Let $C\subset Z$ be the positive-real cycle given by $X_1>0, X_2>0, Y>0$. The cycle $C$ is homeomorphic to the image $\Log_t(C) \subset \R^3$ under the map $\Log_t$: 
\[
\Log_t(C) = \{(x_1,x_2,y) \in \R^3: x_1+x_2 = \log_t(1+t^y) \} 
\] 
The image $\Log_t(C)$ is very close to the union $(x_1+x_2 =\min(0,y))$ of two faces of the tropical amoeba. We now compute the period of $C$ with respect to the holomorphic volume form $\Omega$. The volume form $\Omega$ restricted to $C$ can be identified with the affine area form in the coordinates $(x_i,y)$: 
\[
\Omega|_C = (-\log t)^2 dx_1 \wedge dy = (- \log t)^2 dx_2 \wedge dy. 
\]
To make the computation finite, we consider the integral over the  finite region $R$ delimited by the affine-linear equations (with $a_1,a_2,b$ positive constant):
\[
R = \{(X_1,X_2,Y) \in C : x_1 \le -a_1,\ x_2  \le -a_2, \ -b\le y \le b  \}
\] 
The corresponding polytope in the affine chart is shown in green colour in Figure \ref{fig:focus-focus}. The affine area of this polytope is $2 (a_1+a_2) b - b^2/2$. On the other hand, the actual shape of $R$ in the coordinate $(x_1,y)$ slightly differs from this polytope because of the error in tropicalization, see Figure \ref{fig:actual_region}. We have 
\begin{align*} 
\int_{R\subset C} \Omega & = \int_{\substack{x_1+x_2 = \log_t(1+t^y) \\ x_i \ge -a_i, |y|\le b}} (-\log t)^2 dx_1 dy \\
& = (-\log t)^2 \int_{-b}^b (a_2+\log_t(1+t^y)- (-a_1))dy \\ 
& =  (-\log t)^2 \left( 
\int_{-b}^b (a_1+a_2+\min(0,y)) dy - \int_{-b}^b (-\log_t(1+t^y) + \min(0,y)) dy \right) \\ 
& = (-\log t)^2 (\text{area of the polytope}) - \zeta(2) + O(t^b).  
\end{align*} 
where $\zeta(2)$ arises from the `error in tropicalization' integral (area of the blue region):  
\begin{equation} 
\label{eq:2dim_error}
\zeta(2) = (-\log t)^2 \int_{-\infty}^\infty (-\log_t(1+t^y) + \min(0,y)) dy.  
\end{equation} 
This means that the singularity contributes $- \zeta(2)$ to the period integral. 
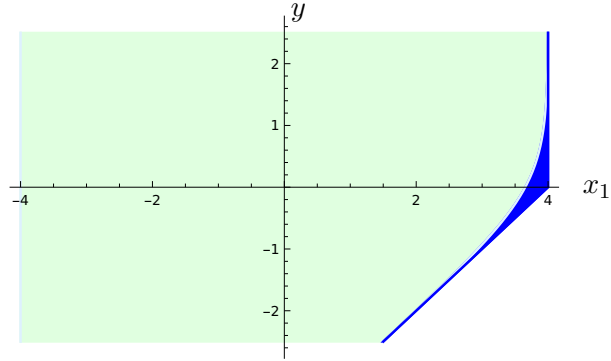
\begin{figure}[ht] 
\centering 
\includegraphics[scale=0.8]{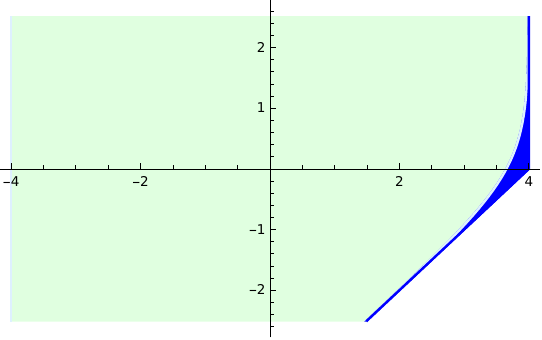}
\begin{picture}(0,0)
\put(10,65){\makebox(0,0){$x_1$}} 
\put(-103,132){\makebox(0,0){$y$}}
\end{picture} 
\caption{The region $R$ in coordinates $(x_1,y)$. It is given by $-a_1\le x_1 \le a_2 + \log_t(1+t^y)$, $|y|\le b$. This differs from the green polytope in Figure \ref{fig:focus-focus} by the blue region.} 
\label{fig:actual_region} 
\end{figure} 
\end{example} 

\begin{example} 
Let $(Z_t,\Omega_t)$ be the family of affine K3 surfaces 
\begin{align*} 
Z_t & = \{(W_1,W_2,W_3) \in (\C^\times)^3: t (W_1+ W_2 + W_3 + \frac{1}{W_1W_2W_3}) = 1\} \\ 
\Omega_t & = \frac{d \log W_1 \wedge d \log W_2 \wedge d \log W_3}{d (t (W_1+W_2+ W_3 +
\frac{1}{W_1W_2W_3}))} 
\end{align*}
This can be compactified to a K3 surface $\overline{Z_t}$ so that $\Omega_t$ extends to a nowhere vanishing holomorphic 2-form. This family is mirror to a quartic K3 surface $Y\subset \PP^3$ equipped with a symplectic form in the class $(-\log t) c_1(\PP^3)$. The tropical amoeba of $Z_t$ is given by the singular locus of $\min(w_1+1, w_2+1, w_3+1, 1-w_1-w_2-w_3,0)$ as in Figure \ref{fig:trop_K3}. Observe that the compact chamber 
\[
\Delta = \{(w_1,w_2,w_3) \in \R^3: w_1\ge -1, w_2\ge -1, w_3 \ge -1, w_1+w_2+w_3 \le 1\}
\]
bounded by the tropical amoeba is the same as the moment polytope of $\PP^3$ with respect to $c_1(\PP^3)$. By collapsing\footnote{See \cite{Yamamoto} for collapsing tropical hypersurfaces.} all the unbounded faces of the tropical amoeba (not contained in $\Delta$), we should get a torus fibration on $\overline{Z_t}$ over the sphere $\partial \Delta \cong S^2$. This fibration has singularities, since there is no uniform direction to collapse the unbounded faces, unlike the case of Example \ref{exa:elliptic}. (In fact, we have many ways to collapse, yielding different fibration structures on $\overline{Z_t}$.)

\begin{figure}[ht]
\centering 
\includegraphics[scale=0.5]{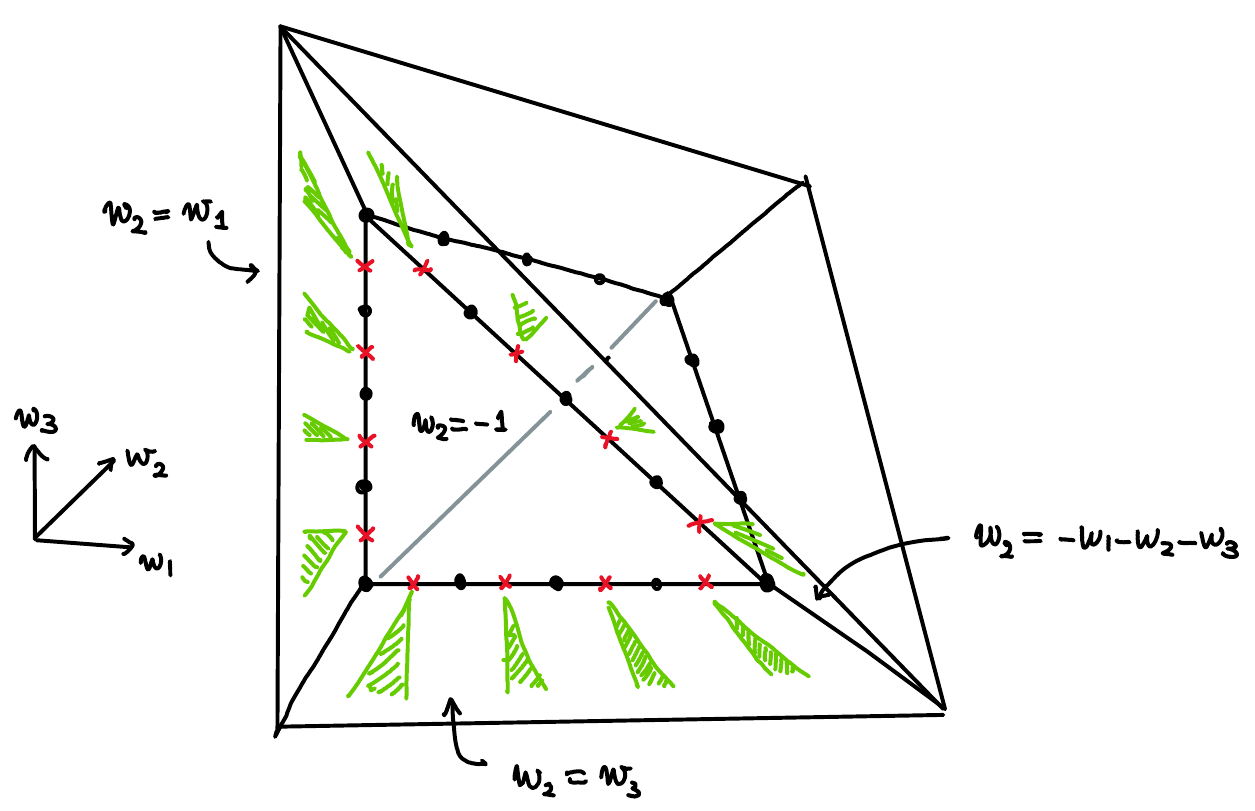}
\caption{Tropical K3 surface. This is composed of the boundary of the 3-simplex $\Delta$ and 6 unbounded faces. The green regions are the $\Log_t$-images of neighbourhoods of singularities (as in Figure \ref{fig:2d_sing}) and the red crosses represent singular points of the affine structure on $\partial \Delta$.}
\label{fig:trop_K3}
\end{figure} 

As before, we consider the period of the positive real cycle $C_t = \{(W_1,W_2,W_3) \in Z_t : W_1 >0, W_2>0, W_3>0\}$. It is homeomorphic to its image under $\Log_t$:
\[
\Log_t(C_t) = \{(w_1,w_2,w_3) \in \R^3: t^{w_1+1} + t^{w_2+1} + t^{w_3+1} + t^{1-w_1-w_2-w_3} = 1\}. 
\]
As $t\to +0$, $\Log_t(C_t)$ rapidly approaches $\partial\Delta$. We compute the period integral locally over the tropical base $\partial \Delta$, reducing the computation to the local model from Example \ref{exa:2d_local}. Since the fibration has 24 singular points (see below), we expect to have the following asymptotics 
\begin{equation}
\label{eq:K3_asymp}
\int_{C_t} \Omega_t = (-\log t)^2 \text{(affine area of $\partial \Delta$)} - 24 \zeta(2) + O(t^\epsilon). 
\end{equation} 
The number $-24 \zeta(2)$ also arises from the $\hGamma$-class of a quartic $Y$: $\int_Y \hGamma_Y = - \zeta(2) \int_Y c_2(Y)$. We shall explain the details below. We introduce a $\Z$-affine structure (with singularities) on $\partial \Delta$ as follows (see \cite{Gross-Siebert}). 
\begin{itemize} 
\item on the interior of 2-dimensional faces of $\partial \Delta$, we consider the subspace affine structure; 
\item on a neighbourhood of a lattice point $v$ on edges (1-dimensional face) of $\partial \Delta$, we consider the affine structure induced from the projection $\R^3 \to \R^3/\langle v\rangle$ and;  
\item we need to have a singularity (shown in red cross) somewhere between adjacent lattice points $v_1,v_2$ on an edge because the affine structures induced by the projections $\R^3 \to \R^3/\langle v_i \rangle$, $i=1,2$ are not compatible with each other; we have 4 singularities on each edge, and thus $24=4\times 6$ singularities in total. 
\end{itemize}  
The reason for this affine structure is as follows. Near the interior of the 2-dimensional face given by $(w_1=-1)$ for example, the cycle $C_t \cong \Log_t(C_t)$ is approximated by the affine-linear subspace 
\[
t^{w_1+1} \approx 1 
\]
since the other terms $t^{w_2+1}, t^{w_3+1}, t^{1-w_1-w_2-w_3}$ are much smaller than $t^{w_1+1}$, and the volume form $\Omega_t|_{C_t}$ is approximated by the affine volume form (mulitplied by $(-\log t)^2$): 
\begin{align*}
\Omega_t|_{C_t} &= \frac{d \log t^{w_1+1} \wedge d \log t^{w_2+1} \wedge d \log t^{w_3+1}}{d (t^{w_1+1}+t^{w_2+1} + t^{w_3+1} + t^{1-w_1-w_2-w_3})} \\
& \approx \frac{(-\log t)^3 dw_1 \wedge dw_2 \wedge dw_3}{d(t^{w_1+1})} \qquad \text{(near the face $(w_1=-1) \cap \Delta$)} \\ 
& = (-\log t)^2 d w_2 \wedge d w_3. 
\end{align*} 
Next, around the interior of the edge given by $(w_1=w_2=-1)$ for example, the cycle $C_t$ is approximated by 
\begin{equation} 
\label{eq:approx_Ct} 
t^{w_1+1}+t^{w_2+1} \approx 1 
\end{equation} 
since the other terms $t^{w_3+1}, t^{1-w_1-w_2-w_3}$ are exponentially small. The volume form $\Omega_t|_{C_t}$ is similarly approximated by the affine volume form on $\R^2/\langle v\rangle$ for a lattice point $v$ on the edge. 
\begin{align*} 
\Omega_t|_{C_t} 
& \approx \frac{(-\log t)^3 dw_1 \wedge dw_2 \wedge dw_3}{d(t^{w_1+1}+t^{w_2+1})} 
= (-\log t)^2 \text{(affine volume form on $\R^3/\langle v\rangle$)}.  
\end{align*} 
In fact, if we complete $v$ to a $\Z$-basis $v=v_1, v_2, v_3$ of $\Z^3$ and write $(a_1,a_2,a_3)$ for the linear coordinates on $\R^3$ dual to $(v_1,v_2,v_3)$, we have $w_i(a_1 v) = -a_1$; thus $\parfrac{w_1}{a_1} = \parfrac{w_2}{a_1}=-1$ and 
\begin{align*} 
\frac{(-\log t)^3 dw_1 \wedge dw_2 \wedge dw_3}{d(t^{w_1+1}+t^{w_2+1})} 
& = \frac{(-\log t)^3 da_1 \wedge da_2 \wedge da_3}{d(t^{w_1+1} + t^{w_2+1})} \\
& = \frac{(-\log t)^3 da_2 \wedge da_3}{\partial(t^{w_1+1}+t^{w_2+1})/\partial a_1} 
= (-\log t)^2 da_2 \wedge da_3.  
\end{align*}  
Note that the affine structure depends on the choice of $v$, but the volume form does not. 
Finally, around the vertex $v=(-1,-1,-1)$ of $\partial \Delta$, $C_t$ is approximated by 
\[
t^{w_1+1}+t^{w_2+1}+t^{w_3+1} \approx 1 
\]
and $\Omega_t|_{C_t}$ is also approximated by the affine volume form on $\R^3/\langle v\rangle$. 

We next examine the cycle $C_t$ near the singularity. Let $v_1,v_2$ be adjacent lattice points on the edge $w_1=w_2=-1$. Choose an integral vector $v_3$ such that $(v_1,v_2,v_3)$ is a $\Z$-basis of $\Z^3$ 
and let $(a_1,a_2,a_3)$ be the linear coordinates on $\R^3$ dual to $(v_1,v_2,v_3)$. We may assume that $v_3$ is parallel to the face $w_1=-1$ (i.e.~$w_1(v_3)=0$) by adding a linear combination of $v_1$ and $v_2$ if necessary. The coordinates $w_1,w_2,w_3$ can then be written as a $\Z$-linear combination of $a_1,a_2,a_3$; examining the values at $v_1,v_2,v_3$ we find that the coordinate change is of the form: 
\begin{align*} 
w_1 & = -a_1 -a_2 \\ 
w_2 & = -a_1 -a_2 + m a_3 \\
w_3 & = j a_1 + k a_2 + l a_3 
\end{align*} 
for some $j,k,l,m\in \Z$. Moreover, since $w_1-w_2$ is a non-zero primitive covector, we have $m=\pm 1$. By flipping the sign of $a_3$ (or equivalently $v_3$), we may assume that $m=1$. Recall from \eqref{eq:approx_Ct} that $C_t$ is given by $t^{w_1+1} + t^{w_2+1} \approx 1$ around the edge. This can be rewritten as: 
\[
1 + t^{a_3} \approx t^{a_1+a_2-1} = t^{a_1-\frac{1}{2}} \cdot t^{a_2 - \frac{1}{2}}.  
\]
Setting $X_1 = t^{a_1-\frac{1}{2}}$, $X_2 = t^{a_2-\frac{1}{2}}$, $Y= t^{a_3}$, this equation can be identified with the local model from Example \ref{exa:2d_local}. The holomorphic volume form $(-\log t)^2 da_1\wedge da_2 = d \log X_1 \wedge d\log Y$ is also the same as in Example \ref{exa:2d_local}. If we introduce a Lagrangian fibration near $\{v_1,v_2\}$ by locally identifying $Z_t$ with the local model in Example \ref{exa:2d_local} by this identification 
(i.e.~$W_1= (tX_1X_2)^{-1}, W_2 = Y (tX_1X_2)^{-1}, W_3 = X_1^jX_2^k Y^l t^{(j+k)/2}$), we see that the singularity appears at $(a_1,a_2,a_3)= (\frac{1}{2},\frac{1}{2}, 0)$ on the base, which corresponds to the mid-point of $v_1$ and $v_2$. At each singular point, we have the error term $-\zeta(2)$ as we saw in Example \ref{exa:2d_local} and arrive at the asymptotics \eqref{eq:K3_asymp}. 

We note that the choice of an (approximate) torus fibration, in particular, the choice of positions of singularities on $\partial \Delta$ has certain arbitrariness in the above discussion. We introduced singularities on $\partial \Delta$ for the purpose of computing periods (and we had to since we cannot cover $\partial \Delta$ by compatible affine charts) and the error terms in tropical approximation occurred from those singularities. 
\end{example} 

\begin{example} 
We have seen that $\zeta(2)$ arises from 2-dimensional singularities through the `error in tropicalization' integral \eqref{eq:2dim_error}. In dimension three, $\zeta(3)$ arises from the following `error in tropicalization' integrals  
(see \cite[Eqn (22), Proposition 4.5]{AGIS}): 
\begin{align*}
& \zeta(3) = \frac{1}{2} (-\log t)^3 \int_{-\infty}^\infty\left( \left(\log_t(1+t^y)\right)^2 - \left(\min(0,y)\right)^2 \right) dy  \\ 
\intertext{and} 
& (-\log t)^3 \int_U (- \log_t(1+t^{y_1} + t^{y_2}) + \min(0,y_1,y_2)) dy_1 dy_2 \\
& \qquad \qquad \quad = (-\log t) \ell \zeta(2)  +  \chi_U(0) \zeta(3)+O(1/(-\log t)) 
\end{align*} 
where $U$ is a bounded domain in $\R^2$ such that $\partial U$ intersects every stratum of the tropical curve $L= \Sing(\min(0,y_1,y_2))$ transversally, $\ell $ is the total affine length of $U\cap L$ and 
\[
\chi_U(0)=\begin{cases} 1 & \text{if $0\in U$,} \\ 
0 & \text{otherwise.} 
\end{cases}
\]
\end{example}

\begin{remark} 
In dimension one, the relationship between periods and affine-length of tropical curves was studied in \cite{Mikhalkin-Zharkov, Iwao}. For toric Calabi-Yau hypersurfaces, the relationship between periods and the radiance obstruction was studied in \cite{Yamamoto}. Ruddat and Siebert \cite{Ruddat-Siebert} computed periods of cycles fibering over 1-dimensional curves on the SYZ base, in the framework of Gross-Siebert program. 
\end{remark}

\providecommand{\arxiv}[1]{\href{http://arxiv.org/abs/#1}{arXiv:#1}}

\end{document}